\documentclass[12pt]{amsart}
\usepackage{amsmath}
\usepackage{amsfonts,amssymb,amsthm, txfonts, pxfonts}
\def\struckint{\mathop{%
\def\mathpalette##1##2{\mathchoice{##1\displaystyle##2}%
  {##1\textstyle##2}{##1\scriptstyle##2}{##1\scriptscriptstyle##2}}%
\mathpalette
{\vbox\bgroup\baselineskip0pt\lineskiplimit-1000pt\lineskip-1000pt
\halign\bgroup\hfill$}
{##$\hfill\cr{\intop}\cr\diagup\cr\egroup\egroup}%
}\limits}

\newcommand{\Nmin}{\mathrm{Nmin\:}}

\newtheorem{theorem}{Theorem}
\newtheorem{lemma}[theorem]{Lemma}

\newtheorem{corollary}[theorem]{Corollary}

\newtheorem{definition}[theorem]{Definition}

\theoremstyle{remark}

\newtheorem{remark}[theorem]{Remark}

\newtheorem{example}[theorem]{Example}

\newcommand{\essn}{\mathbb{S}^n}

\newcommand{\reals}{\mathbb{R}}

\providecommand{\mindex}[1]{\mathbf{#1}}
\providecommand{\mindexc}[1]{\overline{\mindex{#1}}}
\DeclareMathOperator{\vol}{vol}
\DeclareMathOperator{\Hom}{Hom}
\DeclareMathOperator{\rank}{rank}
\DeclareMathOperator{\diag}{diag}

\begin{document}


\title{Surface area and other measures of Ellipsoids}

\author{Igor Rivin}


\address{Department of Mathematics, Temple University, Philadelphia}

\email{rivin@math.temple.edu}

\thanks{The author is supported by the NSF DMS. Parts of this paper appeared
  as the preprint \cite{riv1}; the author would like to thank Warren~D.~Smith
  on comments on a previous version of this paper.The author would also like
  to thank Princeton University, New York University, and Unversit\'e Paul
  Sabatier for their hospitality.The author would also like
  to thank Franck Barthe for simplifying the arguments in section
  \ref{fbounds},  and thus obtaining the sharp bounds presented in that
  section, and also in Theorem \ref{theestimate} (which uses the same
  argument). Barthe also found an integral-geometric proof of much of Theorem 
  \ref{measproj} } 

\date{\today}

\keywords{ellipsoid, surface area, expectation, integration, large dimension,
  Lindberg conditions, Lauricella hypergeometric function, harmonic mean}

\subjclass{52A38; 60F99; 58C35}

\begin{abstract}
We begin by studying the surface area of an ellipsoid in $\mathbb{E}^n$ as the
function of 
the lengths of the semi-axes.  We write down an explicit formula as an
integral over $\mathbb{S}^{n-1},$ use this formula to derive convexity
properties of the surface area, to give sharp estimates
for the surface area of a large-dimensional ellipsoid,
to produce asymptotic formulas for the surface area and the
\emph{isoperimetric ratio} of an ellipsoid in large dimensions, and to give an
expression for the surface in terms of the Lauricella hypergeometric
function.
We then write down general formulas for the volumes of projections of
ellipsoids, and use them to extend the above-mentioned results to give
explicit and approximate formulas for the higher integral mean curvatures of
ellipsoids. 
Some of our results can be expressed as \emph{isoperimetric} results for
higher mean curvatures.
\end{abstract}

\maketitle

\section*{Introduction}
We study the \emph{mean curvature integrals} of
an ellipsoid $E$ in $\mathbb{E}^n$  as functions of the lengths of its
\emph{semiaxes} -- the $0$-th mean curvature integral is simply the surface
area of the ellipsoid $E.$ The goal is to study the properties of these
integral mean curvatures as functions of the lengths of the (semi)axes of the
ellipsoid. We derive explicit formulas, very good approximations, and
asymptotic results. In addition, some of our results can be viewed as
isoperimetric results for ellipsoids, and we conjecture generalizations to
hold for arbitrary convex bodies.

In detail,  we first write down a formula (\eqref{iratio}) expressing the surface
area of $E$ in terms of an integral of a simple function over the sphere
$\mathbb{S}^{n-1}.$ This formula will be used to deduce a number of results:
\begin{enumerate}
\item The ratio of the surface area to the volume of $E$ (call
this ratio
  $\mathcal{R}(E)$) is a \emph{norm} on
  the vectors of inverse semi-axes. (Theorem \ref{isnorm}).
\item By a simple transformation (introduced for this purpose in \cite{riv1},
  though doubtlessly known for quite some time) $\mathcal{R}(E)$ can be
  expressed as a moment of a sum of independent Gaussian random
  variables; this transformation can be used to evaluate or estimate quite a
  number of related spherical integrals (see Section \ref{sphints}).
\item Sharp bounds (\eqref{sharpineq}) on the ratio of
  $\mathcal{R}(E)$ to the $L^2$   norm of the vectors of inverses of
  semi-axes are derived.
\item We write down a very simple  asymptotic formula
  (Theorem \ref{asymp}) for the surface area of an ellipsoid of a very large
  dimension   with ``not too different'' axes. In particular, the formula
  holds if the ratio of the lengths of any two semiaxes is bounded by some
  fixed constant (Corollary \ref{asymp1}).
\item Finally, we give an identity expressing the surface area of $E$  as
a linear combination of Lauricella hypergeometric functions.
\end{enumerate}
We then go on to give similar explicit formulas for the higher mean curvature
integrals of ellipsoids, by first computing the volumes of projections of
ellipsoids onto subspaces (Sections \ref{extalg} and \ref{howtovol}) and then
writing down a simple approximation (Theorem \ref{theestimate}) for the $k$-th
integral mean curvature of an ellipsoid. This estimate (for a fixed $k$) does
not differ from the true value of the $k$-th mean curvature by more than a
(dimension independent) constant factor. The worst possible functional
dependence of our estimate on the dimension $n$ is $O(n^{1/4}),$ which comes to
pass when $k = n/2.$ Our estimates on the error are sharp. Unfortunately, it
seems difficult to derive a law of large numbers (as we describe above for the
surface area -- the $0$-th integral mean curvature).

In Section \ref{history} we comment on some historical antecedents of
our work, and in Section \ref{isoperimetric} we interpret some of our
inequalities as \emph{isoperimetric} inequalities.

\medskip\noindent
\textbf{Notation.} Let $(S, \mu)$ be a measure space with $\mu(S) < \infty.$
We will use the notation
\[
\fint_S f(x) \,d\mu \stackrel{\text{def}}{=} \dfrac{1}{\mu(S)} \int_S
f(x) d\mu.
\]
In addition, we shall denote the area of the unit sphere $\mathbb{S}^n$ by
$\omega_n$ and we shall denote the volume of the unit ball $\mathbb{B}^n$ by
$\kappa_n.$

\section{Cauchy's formula}
Let $K$ be a convex body in $\mathbb{E}^n.$ Let $u \in \mathbb{S}^{n-1}$ be
a unit vector, and let us define $V_u(K)$ to be the (unsigned)
$n-1$-dimensional volume of the orthogonal projection of $K$ in the direction
$U.$ \emph{Cauchy's formula} (see \cite[Chapter 13]{santalo}) then states that
\begin{equation}
\label{cauchy}
\boxed{V_{n-1}(\partial K) = \dfrac{n-1}{\omega_{n-2}}
\int_{\mathbb{S}^{n-1}}V_u(K)\,d\sigma = \\
(n-1)\dfrac{\omega_{n-1}}{\omega_{n-2}} \fint_{\mathbb{S}^{n-1}}
  V_u(K)\,d\sigma,}
\end{equation}
where $d\sigma$ denotes the standard area element on the unit sphere.

In the case where $K=E$ is an ellipsoid, given by
\[
E = \{\mathbf{x}\in \mathbb{E}^n \quad | \quad \sum_{i=1}^n q_i^2 x_i^2 = \leq
1\}
\]
the volume $V_u$ of the projection is computed in Example \ref{dimn1} as a
special case of more general projection results:
\begin{equation}
\label{conform}
V_u(E) = \kappa_{n-1} \dfrac{\sqrt{\left(\sum_{i=1}^n u_i^2
    q_i^2\right)}}{\prod_{i=1}^n q_i}.
\end{equation}

Since
\[
V_n(E) = \dfrac{\kappa_n}{\prod_{i=1}^n q_i},
\]
we can rewrite Cauchy's formula \eqref{cauchy} for $E$ in the form:
\begin{equation}
\label{iratio}
\boxed{\mathbb{R}(E) \stackrel{\text{def}}{=} \dfrac{V_{n-1}(\partial
    E)}{V_n(E)} =  n \fint_{\mathbb{S}^{n-1}} \sqrt{\sum_{i=1}^n u_i^2
    q_i^2}\,d\sigma,}
\end{equation}
where $\mathbb{R}(E)$ is the \emph{isoperimetric ratio} of $E.$
\begin{theorem}
\label{isnorm}
The ratio $\mathbb{R}(E)$ is a \emph{norm} on the vectors $q$ of lengths of
semiaxes ($q = (q_1, \dots, q_n).$)
\end{theorem}
\begin{proof}
The integrand in the formula \eqref{iratio} is a norm.
\end{proof}
\begin{corollary}
\label{pineq}
There exist constants $c_{n, p}, C_{n,p},$ such that
\[
c_{n, p} \| q \|_p \leq \mathbb{R}(q) \leq C_{n, p} \|q\|_p,
\]
where $\|q\|_p$ is the $L^p$ norm of $q.$
\end{corollary}
\begin{proof} Immediate (since all norms on a finite-dimensional Banach space
  are equivalent).
\end{proof}
In the sequel we will find sharp bounds on the constants $c_{n, p}$
and $C_{n, p},$ but for the moment observe that if $a_i = 1/q_i,\quad i=1,
\dots, n,$ then
\begin{equation}
\label{reverse}
\|q\|_p = \prod_{i=1}^n q_i \sigma_{n-1}^{1/p}(a_1^p, \dots,
a_n^p),
\end{equation}
where $\sigma_{n-1}$ is the $n-1$-st elementary symmetric function. In
particular, for $p=1,$ Corollary \ref{pineq} together with Eq. \eqref{reverse}
gives the estimate of \cite{riv2} (only for the $0$-th mean curvature integral
and with (for now) ineffective constants -- the latter part will be remedied
directly).
To exploit the formula \eqref{iratio} fully, we will need a digression on
computing spherical integrals.

\section{Spherical integrals}
\label{sphints}
In this section we will prove the following easy but very useful Theorem:
\begin{theorem}
\label{homothm}
Let $f(x_1, \dots, x_n)$ be a homogeneous function on $\mathbb{E}^n$ of degree
$d$ (in other words, $f(\lambda x_1, \dots, \lambda x_n) = \lambda^d f(x_1,
\dots, x_n)$.) Then
\[
\Gamma\left(\frac{n+d}2\right)
\fint_{\mathbb{S}^{n-1}} f d \sigma =
  \Gamma\left(\frac{n}{2}\right)\mathbb{E}\left(f(\mathbf{X}_1, \dots,
  \mathbf{X}_n)\right),\]
where $\mathbf{X}_1, \dots, \mathbf{X}_n$ are independent random variables
with probability density $e^{-x^2}.$
\end{theorem}
\begin{proof}
Let
\[E(f) = \fint_{\mathbb{S}^{n-1}} f(x)\,d\sigma,\]
and let $N(f)$ be defined as $\mathbb{E}(f(\mathbb{X}_1, \dots,
\mathbb{X_n})),$
where $\mathbb{X}_i$ is a Gaussian random variable with mean $0$ and variance
$1/2,$ (so with probability density $\mathfrak{n}(x) = e^{-x^2},$) and
$\mathbb{X}_1, \dots, \mathbb{X}_n$ are independent.
By definition,
\begin{equation}
\label{ndef}
N(f)(n) = c_n \int_{\mathbb{E}^n} \exp\left(-\sum_{i=1}^n x_i^2\right)
f(x_1, \dots, x_n) \,d x_1 \dots d x_n,
\end{equation}
where $c_n$ is such that
\begin{equation}
\label{norm}
c_n \int_{\mathbb{E}^n} \exp\left(-\sum_{i=1}^n x_i^2\right)\, d x_1 \dots d
x_n = 1.
\end{equation}
We can rewrite the expression \eqref{ndef} for $N(f)$ in polar coordinates as
follows (using the homogeneity of $f$):
\begin{equation}
\Nmin(n) = c_n \vol \mathbb{S}^{n-1}\int_0^\infty e^{-r^2} r^{n+d-1} E(f)
\,d r =  c_n E(f) \int_0^\infty e^{-r^2} r^{n+d-1}\, dr.
\end{equation}
 Since, by the substitution $u = r^2,$
\[\int_0^\infty e^{-r^2} r^{n+d-1}\, dr = \frac{1}{2}\int_0^\infty e^{-u}
u^{(n+d-2)/2} d u = \frac{1}{2}\Gamma\left(\frac{n+d}2\right).\]
and Eq. (\ref{norm}) can be rewritten in polar coordinates as
\[
1=c_n \vol \mathbb{S}^{n-1} \int_0^\infty r^{n-1} dr = \dfrac{c_n \vol
  \mathbb{S}^{n-1}}2 \Gamma\left(\frac{n}2\right),\]  we see that
\[\Gamma\left(\frac{n+d}2\right)E(f) = \Gamma\left(\frac{n}{2}\right)N(f).\]
\end{proof}
\begin{remark}
\label{gammarat}
In the sequel we will frequently be concerned with asymptotic results, so it
is useful to state the following asymptotic formula (which follows immediately
from Stirling's formula):
\begin{equation}
\label{gammaeq}
\lim_{x\rightarrow \infty} \dfrac{\Gamma(x+y)}{\Gamma(x)(x+y)^y} = 1.
\end{equation}
It follows that for large $n$ and fixed $d,$
\begin{equation}
\label{factorasymp}
\dfrac{\Gamma\left(\frac{n}2\right)}{\Gamma\left(\frac{n+d}2\right)} \sim
\left(\frac2{n+d}\right)^{d/2}.
\end{equation}
\end{remark}
\section{An explicit formula for the surface area}
The Theorem in the preceding section can be used to give explicit formulas
for the surface area of an ellipsoid (this formula will not be used in the
sequel, however). Specifically, in the book \cite{quad} there are formulas for
the moments of of random variables which are quadratic forms in Gaussian
random variables. We know that for our ellipsoid $E,$
\[\mathbb{R}(E) = n \fint_{\mathbb{S}^{n-1}} \sqrt{\sum_{i=1}^n u_i^2
  q_i^2}\,d\sigma = n
  \dfrac{\Gamma\left(\frac{n}2\right)}{\Gamma\left(\frac{n+1}2\right)}
\mathbb{E}\left(\sqrt{q_1^2 \mathbb{X}_1 +  \dotsb + q_n^2
  \mathbb{X}_n}\right),\] 
where $\mathbb{X}_i$ is a Gaussian with variance
  $1/2.$ The expectation in the last expression is the $1/2$-th moment of the
  quadratic form in Gaussian random variables, and so the results of
  \cite[p.~62]{quad} apply verbatim, so that we obtain:
\begin{equation}
\label{quadform1}
\mathbb{R}(E) = n
\dfrac{\Gamma\left(\frac{n}2\right)}{\Gamma\left(\frac{n+1}2\right)\Gamma\left(\frac12\right)}\sqrt{\alpha}
\int_0^\infty \dfrac{1}{\sqrt{z}} \sum_{j=1}^n \dfrac{q_j^2}{2(1+\alpha z
  q_j^2)} \left(\prod_{j=1}^n(1-q_j^2 z)\right)^{-1/2}\,dz;
\end{equation}
note that $\alpha$ in the above formula can be any positive number (as long as
 $|1-\alpha q_j^2| < 1,$ for all $j.$

This can also be expressed in terms of special functions.
 First, we need a definition:
\begin{definition}
Let $a, b_1, \dots, b_n, c, x_1, \dots, x_n$ be complex numbers, with
$|x_i| < 1,\quad i=1, \dots, n$, $\Re a > 0,$ $\Re(c-a) > 0.$ We then define
the \emph{Lauricella Hypergeometric Function} $F_D(a; b_1, \dots, b_n; c; x_1,
\dots, x_n)$ as follows:
\begin{multline}
F_D(a; b_1, \dots, b_n; c; x_1, \dots, x_n) = \\
\dfrac{\Gamma(c)}{\Gamma(a)\Gamma(c-a)} \int_0^1
u^{a-1}(1-u)^{c-a-1}\prod_{i=1}^n(1-ux_i)^{-b_i}\,d u.
\end{multline}
We also have the series expansion:
\begin{multline}F_D(a; b_1, \dots, b_n; c; x_1, \dots, x_n) = \\
\sum_{m_1 = 0}^\infty \cdots \sum_{m_n = 0}^\infty
\dfrac{(a)_{m_1 + \cdots + m_n} \prod_{i=1}^n (b_i)_{m_i}}{(c)_{m_1 + \cdots +
    m_n}} \prod_{i=1}^n \dfrac{x_i^{m_i}}{m_i!},
\end{multline}
valid whenever $|x_i| < 1, \forall i.$
\end{definition}
Now, we can write
\begin{multline}
\label{quadform2}
\mathbb{R}(E) = n
\dfrac{\Gamma^2\left(\frac{n}2\right)}{\Gamma^2\left(\frac{n+1}2\right)}
\sqrt{\alpha}\times\\
 \sum_{j=1}^n \frac{q_j^2}{2}
F_D\left(1/2; \eta_{1j}, \dotsc, \eta_{nj}; \frac{n+1}2; 1-\alpha q_1^2, \dotsc,
1-\alpha q_n^2\right),
\end{multline}
where $\eta_{ij} = 1/2 + \delta_{ij},$ and $\alpha$ is a positive parameter
satisfying $|1-\alpha q_j^2| < 1.$

\section{Laws of large numbers}
Many of the results in this section will require the following basic lemmas.
\begin{lemma}
\label{twomo}
Let $F_1, \dotsc, F_n, \dotsc$ be a sequence of probability distributions
whose first moments converge to $\mu$ and whose second moments converge to
$0.$ then $F_i$ converge to the Dirac delta function distribution centered on
$\mu.$
\end{lemma}
\begin{proof}
Follows immediately from Chebyshev's inequality.
\end{proof}
\begin{lemma}
\label{monotone}
Suppose the distributions $F_1, \dotsc, F_n, \dotsc$ converge to the
distribution $F,$ and the expectations of $|x|^\alpha$ with respect to $F_1,
\dots, F_n, \dotsc$ are bounded. Then the expectation of $|x|^\beta,$ $0
\leq \beta < \alpha$ converges to the expectation of $|x|^\beta$ with respect
to $F.$
\end{lemma}
\begin{proof}
See \cite[pp.~251-252]{feller2}.
\end{proof}
\begin{theorem}
\label{basiclaw}
Let $\mathbb{Y}_1, \dotsc, \mathbb{Y}_n, \dotsc$ be \emph{independent} random
variables with means
$0<\mu_1, \dotsc, \mu_n, \dotsc < \infty$ and variances $\sigma_1^2, \dotsc,
\sigma_n^2, \dotsc < \infty$ such that
\begin{equation}
\label{lovar}
\lim_{n \rightarrow \infty} \dfrac{\sum_{i=1}^n \sigma_i^2}{\left(\sum_{i=1}^n
  \mu_i \right)^2} = 0.
\end{equation}
Then
\[
\lim_{n\rightarrow \infty}\mathbb{E}\left(
\dfrac{\mathbb{Y}_1 + \dotsb + \mathbb{Y}_n}{\sum_{i=1}^n \mu_i}\right)^\alpha
= 1,\] for  $\alpha < 2.$
\end{theorem}
\begin{proof}
Consider the variable
\[
\mathbb{Z}_n = \dfrac{\sum_{i=1}^n \mathbb{Y}_i}{\sum_{i=1}^n \mu_i}.
\]
It is not hard to compute that
\[
\sigma^2(\mathbb{Z}_n) = \dfrac{\sum_{i=1}^n \sigma_i^2}{\left(\sum_{i=1}^n
  \mu_i \right)^2},
\]
while
\[\mu(\mathbb{Z}_n) = 1,\]
so by assumption \eqref{lovar} and Lemma \ref{twomo} $\mathbb{Z}_n$ converges
in distribution to the delta function centered at $1.$ The conclusion of the
Theorem then follows from Lemma \ref{monotone}.
\end{proof}
\begin{lemma}
\label{expcomp}
Let $\mathbb{X}$ be normal with mean $0$ and variance $1/2$ (so probability
density $e^{-x^2}/\sqrt{\pi}.$) Then
\[\mathbb{E}(|\mathbb{X}|^p) =
\dfrac{\Gamma\left(\frac{p+1}2\right)}{\sqrt{\pi}}.\]
\end{lemma}
\begin{proof}
\[\mathbb{E}(|\mathbb{X}|^p) = \frac{2}{\sqrt{\pi}}\int_0^\infty x^p
e^{-x^2}\,dx =
\frac{1}{\sqrt{\pi}}\int_0^\infty u^{(p-1)/2} e^{-u} =
\dfrac{\Gamma\left(\frac{p+1}2\right)}{\sqrt{\pi}}.\]
\end{proof}
\begin{theorem}
\label{normasymp}
\[
\fint_{\mathbb{S}^{n-1}}\|u\|_p \,d \sigma \sim
\dfrac{\Gamma\left(\frac{n}2\right)}{\Gamma\left(\frac{n+1}2\right)} \left(n
  \dfrac{\Gamma\left(\frac{p+1}2\right)}{\sqrt{\pi}}\right)^{\frac{1}{p}}.
\]
\end{theorem}
\begin{proof}
This follows immediately from the $1$-homogeneity of the $L^p$ norm, the
results of Section \ref{sphints}, Theorem \ref{basiclaw}, and Lemma
\ref{expcomp}.
\end{proof}
\subsection{Asymptotics of $\mathbb{R}(E).$}
\begin{theorem}
\label{asymp}
Let $q_1, \dotsc, q_n, \dotsc$ be a sequence of positive numbers such that
\[\lim_{n\rightarrow \infty} \dfrac{\sum_{i=1}^n q_i^4}{\left(\sum_{i=1}^n
  q_i^2 \right)^2} = 0.
\]
Let $E_n$ be the ellipsoid in $\mathbb{E}^n$ with semiaxes $a_1 = 1/q_1,
\dotsc, a_n = 1/q_n.$ Then
\[
\lim_{n\rightarrow \infty}
  \dfrac{\Gamma\left(\frac{n+1}2\right)}{\Gamma\left(\frac{n}2\right)}
\dfrac{\mathbb{R}(E_n)}{n\sqrt{\frac12\sum_{i=1}^n q_i^2}} = 1.
\]
\end{theorem}

\begin{proof}
The Theorem follows immediately from Theorem \ref{basiclaw} and the results of
Section \ref{sphints}.
\end{proof}
\begin{corollary}
\label{asymp1}
Let $a_1, \dotsc, a_n, \dotsc$ be such that $0< c_1 \leq a_i/a_j \leq c_2 <
\infty,$ for any $i, j.$ Let $E_n$ be the ellipsoid with major semi-axes
$a_1, \dotsc, a_n.$ Then
\[
\lim_{n\rightarrow \infty}
  \dfrac{\Gamma\left(\frac{n+1}2\right)}{\Gamma\left(\frac{n}2\right)}
\dfrac{\mathbb{R}(E_n)}{n\sqrt{\frac12\sum_{i=1}^n \frac{1}{a_i^2}}} = 1.
\]
\end{corollary}
\begin{proof}
The quantities $q_1 = 1/a_1, \dotsc, q_n = 1/a_n, \dotsc$ clearly satisfy the
hypotheses of Theorem \ref{asymp}
\end{proof}

\section{General bounds on $\mathbb{R}(E)$}
\label{fbounds}
We know that $\mathbb{R}(E)$ is a norm on the vector
 $\mathbf{q}=(q_1, \dotsc,q_n)$  -- let us agree to write
\[
\|\mathbf{q}\|_\mathbb{R} \stackrel{\text{def}}{=}
\dfrac{\mathbb{R}(E)}{n}  = \fint_{\mathbb{S}^{n-1}}
\sqrt{\sum_{i=1}^n q_i^2 x_i^2} \,d \sigma.
\]
where $\mathbf{q}$ is the vector of inverses of the major semi-axes of $E.$

We know that
\[ 
c_n \| q\| \leq \|q\|_\mathbb{R} \leq C_n \|q\|,
\] 
for some
dimensional constants $c_n, C_n.$ In this section we
will give sharp  estimates on the constants $c_n$ and $C_n.$

These estimates will depend on the following observation:

\begin{lemma}
\label{convexity} 
Let $\alpha_1, \dots, \alpha_n$ be nonnegative real numbers, and let
\[
f(\alpha_1, \dots, \alpha_n) = \fint_{\mathbb{S}^{n-1}}
 \sqrt{\sum_{i=1}^n \alpha_i x_i^2} d\sigma.
\]
Then $f(\mathbf{\alpha})$ is a concave function of the vector
$\mathbf{\alpha}=(\alpha_1, \dots, \alpha_n).$
\end{lemma}
\begin{proof}
The integrand is concave, since the square root is a concave function. The
integral is thus also concave, as a sum of concave functions.
\end{proof}

\begin{lemma}
\label{sharplem}
The ratio
\[ \dfrac{\|q\|_{\mathbb{R}}}{\|q\|}\] is \emph{maximized} when all of the
$q_i$ are equal; it is \emph{minimized} when $q_2 = \dots = q_n = 0.$
\end{lemma}
\begin{proof}
By homogeneity, we can assume that $\|q\| = 1.$
Now, let $\alpha_i = q_i^2.$ Letting 
\[f(\mathbf{\alpha}) =
\|q\|_{\mathbb{R}},\]
 we see that $f$ is a \emph{symmetric} function, while Lemma
\ref{convexity} tells us that $f(\mathbf{\alpha})$ 
is a \emph{concave} function.  Since the set 
\[ 
S: \sum_{i=1}^n \alpha_i = 1
\] 
is convex, we know that the maximum of $f$ is attained at the point of maximum
symmetry ($\alpha_i = 1/n,$ for all $i$), and the minimum at an extreme point
of $S$ -- by symmetry any extreme point will do, for example $(1, 0, \dotsc,
0).$  
\end{proof}
\begin{corollary}
The minimal value (previously denoted by $c_n$)
 of $\|q\|_{\mathbb{R}}/\|q\|$ equals 
\[
\fint_{\essn} |x_1| d\sigma =
\dfrac{\Gamma\left(\dfrac{n}{2}\right)}{\sqrt{\pi}\Gamma\left(\dfrac{n+1}{2}\right)},    
\]
while the maximal value ($C_n$) of $\|q\|_{\mathbb{R}}/\|q\|$ equals
$1/\sqrt{n},$ or in other words,
\begin{equation}
\label{sharpineq}
\dfrac{\Gamma\left(\dfrac{n}{2}\right)}{\sqrt{\pi}\Gamma\left(\dfrac{n+1}{2}\right)}
\leq
\dfrac{\|q\|_{\mathbb{R}}}{\|q\|} 
\leq \dfrac1{\sqrt{n}}. 
\end{equation}
\end{corollary}

\begin{remark}
The left hand side of equation \eqref{sharpineq} is asymptotic to 
\[
\sqrt{\dfrac2{(n+1)\pi}},\] (by \eqref{factorasymp}) so the ratio $C_n/c_n$
  approaches $\sqrt{2/\pi} 
  = 0.797$ as $n$ goes to infinity.
\end{remark}
\section{Some historical remarks}
\label{history}
The perimeter of an ellipse has been studied since at least Fagnano
(1716). The best approximation has been obtained by Ramanujan in
1914. The surface area of an ellipsoid was studied by Monge
\cite{monge} and
Legendre \cite{legendre} by means of elliptic integrals. Monge also
gave an approximate formula (as a series) which converges as long as
the ellipsoid is not too round; Legendre gave a generally convergent
series. Interesting estimates (also in dimension $3$), especially for
the \emph{mean curvature} of the ellipsoid  were given by
G. Polya and G. Szeg\"o in \cite{posz}. Almost none of the methods in
the references cited above seem to extend to dimension higher than
three. It would be interesting to extend the results and method of the
current article to higher integral mean curvatures of ellipsoids, as
studied in \cite{riv2}.
\section{Higher mean curvatures}
In the sequel we will denote the surface area of the $n$-dimensional unit sphere 
by $\omega_n,$ and the volume of the $n$-dimensional unit ball by $\kappa_n.$
We recall that:

\begin{gather*}
\omega_n = \dfrac{2\pi^{(n+1)/2}}{\Gamma\left(\frac{n+1}2\right)},\\
\kappa_n = \dfrac{\omega_{n-1}}{n} =
\dfrac{2\pi^{n/2}}{\Gamma\left(\frac{n}2\right)}.
\end{gather*}

We will use $M_k(K)$ to denote the integral $k$-th mean curvature of the
boundary of a convex body $K.$ Recall that $M_k(B^n(1)) = \omega_{n-1}.$
($B^n(R)$ is the unit ball of radius $R$ in $\mathbb{E}^n.$
The following result can be found in Santal\'o's book \cite{santalo}: 
\begin{equation}
\label{meancurv}
M_k^{(n)}\left(B^{n-k-1}\right) = \dfrac{\omega_k
  \omega_{n-k-2}}{(n-k-1)\binom{n-1}{k}},
\end{equation}
where $M_k^{(n)}(K)$ denotes the $k$-th integral mean curvature of $K$ viewed
as a convex body in $\mathbb{E}^n.$

We will also need the following:
\begin{theorem}
\label{meancrat}
Let $n-1 k > 1.$ Then
\[
\dfrac{M_k\left(B^n(R)\right)}{M_k^{(n)}\left(B^{n-k-1}(R)\right)} = 
2(k-1) \pi^{3/2}
\dfrac{\Gamma\left(\frac{n+1}2\right)}{\Gamma\left(\frac{k}2\right)\Gamma\left(\frac{n-k}2\right)}.
\]
\end{theorem}

The proof will rely on \emph{Legendre's duplication formula}:
\begin{equation}
\label{legendre}
\Gamma(2 z) = \dfrac{2^{2 z -1}}{\sqrt{\pi}} \Gamma(z) \Gamma(z+1/2).
\end{equation}

\begin{proof}[Proof of Theorem \ref{meancrat}]
Using the formulas for the surface area of the sphere and the $k$-th mean
curvature of $B^k$ we write:
\begin{equation}
\label{calcul1}
\begin{split}
\dfrac{M_k\left(B^n(R)\right)}{M_k^{(n)}\left(B^{n-k-1}(R)\right)} = \\
(n-k-1) \binom{n-1}{k} \dfrac{\omega_{n-1}}{\omega_k\omega_{n-k-2}} = \\
\frac12(n-k-1)\binom{n-1}{k} \pi\dfrac{\Gamma\left(\frac{k+1}2\right)
  \Gamma\left(\frac{n-k-1}2\right)}{\Gamma\left(\frac{n}2\right)} = \\
\frac12(n-k-1) \pi\dfrac{\Gamma\left(\frac{k+1}2\right)
  \Gamma\left(\frac{n-k-1}2\right)}{\Gamma\left(\frac{n}2\right)}
\dfrac{\Gamma(n)}{\Gamma(k-1)\Gamma(n-k)} = \\
\dfrac{\pi}2 \dfrac{\Gamma(n)}{\Gamma\left(\frac{n}2\right)}
\dfrac{\Gamma\left(\frac{k+1}2\right)}{\Gamma(k-1)}
\dfrac{\Gamma\left(\frac{n-k-1}2\right)}{\Gamma(n-k-1)} = \\
\dfrac{(k-1)\pi }4 
\dfrac{\Gamma(n)}{\Gamma\left(\frac{n}2\right)}
\dfrac{\Gamma\left(\frac{k-1}2\right)}{\Gamma(k-1)}
\dfrac{\Gamma\left(\frac{n-k-1}2\right)}{\Gamma(n-k-1)}  
\end{split}
\end{equation}
We have used the assumption that $k > 1$ to factor out $k-1$ in the last line,
and that $n-1>k$ to factor out the $(n-k-1)$ in the first line.

Now we apply Legendre's duplication formula \eqref{legendre} with $z = n/2,$
  $z = (k-1)/2,$ and $z = (n-k-1)/2,$ to get:
\begin{gather*}
\dfrac{\Gamma(n)}{\Gamma\left(\frac{n}2\right)} = 
\dfrac{2^{n-1}}{\sqrt{\pi}} \Gamma\left(\frac{n+1}2\right),\\
\dfrac{\Gamma\left(\frac{k-1}2\right)}{\Gamma(k-1)} = 
\dfrac{\sqrt{\pi}}{2^{k-2}} \dfrac{1}{\Gamma\left(\frac{k}2\right)},\\
\dfrac{\Gamma\left(\frac{n-k-1}2\right)}{\Gamma(n-k-1)} = 
\dfrac{\sqrt{\pi}}{2^{n-k-2}} \dfrac{1}{\Gamma\left(\frac{n-k}2\right)}.
\end{gather*}
Substituting back into \eqref{calcul1}, we see that:
\[
\dfrac{M_k\left(B^n(R)\right)}{M_k^{(n)}\left(B^{n-k-1}(R)\right)} = 
2(k-1) \pi^{3/2}
\dfrac{\Gamma\left(\frac{n+1}2\right)}{\Gamma\left(\frac{k}2\right)\Gamma\left(\frac{n-k}2\right)}.
\]
\end{proof}

\begin{theorem}
\label{meancrat2}
Let $k > 1.$ Then
\begin{equation}
\label{meanceq}
\dfrac{\dfrac{M_k\left(B^n(R)\right)}{M_k^{(n)}\left(B^{n-k-1}(R)\right)}}{\sqrt{\binom{n}{k+1}}}= \pi^{5/4}\dfrac{k-1}{\sqrt{k(k+1)}}
\sqrt{\dfrac{\Gamma\left(\frac{n+1}2\right)}{\Gamma\left(\frac{n}2+1\right)}}
\sqrt{\dfrac{\Gamma\left(\frac{k+1}2\right)}{\Gamma\left(\frac{k}2\right)}}
\sqrt{\dfrac{\Gamma\left(\frac{n-k+1}2\right)}{\Gamma\left(\frac{n-k}2\right)}}
.
\end{equation}
\end{theorem}
\begin{proof}
First, we write 
\[
\binom{n}{k+1} = \dfrac{\Gamma(n+1)}{\Gamma(k+2)\Gamma(n-k)}.
\]
and then, using Legendre's duplication formula and the functional equation of
the $\Gamma$ function:
\begin{gather*}
\Gamma(n+1) = \dfrac{2^n}{\sqrt{\pi}} \Gamma\left(\frac{n+1}2\right)
\Gamma\left(\frac{n}2+1\right),\\
\Gamma(k+2) = (k+1)k\Gamma(k) = (k+1)k\dfrac{2^{k-1}}{\sqrt{\pi}} \Gamma\left(\frac{k+1}2\right)
\Gamma\left(\frac{k}2\right),\\
\Gamma(n-k) = \dfrac{2^{n-k-1}}{\sqrt{\pi}} \Gamma\left(\frac{n-k}2\right)
\Gamma\left(\frac{n-k+1}2\right).
\end{gather*}
The result follows by combining the above with the result of Theorem
\ref{meancrat}.
\end{proof}
\begin{remark}
\label{quotest}
It is not hard to see that for a \emph{fixed} $k,$ the right hand side of
Eq.~\eqref{meanceq} approaches 
\[
C(k) = \pi^{5/4}\dfrac{k-1}{\sqrt{k(k+1)}} \sqrt{\dfrac{\Gamma((k+1)/2)}{\Gamma(k/2)}}
\]
as $n \rightarrow \infty.$ For a \emph{fixed} $m = n-k$ (but both $k$ and $n$
tending to $\infty$), the right hand side approaches 
\[
D(m)=\pi^{5/4}\left(\dfrac{n-k+1}{2}\right)^{1/4}.
\]
Finally, if $n,$ $k,$ and $n-k$ all approach infinity, the expression is
asymptotic to 
\[
B(n, k) = \pi^{5/4}\left(\dfrac{(k+1)(n-k+1)}{2(n+2)}\right)^{1/4}.
\]
It is not hard to see that for a given $n$ $B(n, k)$ is maximized when $k =
n/2,$ in which case 
\[
B(n, n/2) = \left(\dfrac{n+2}8\right)^{1/4}.
\]
\end{remark}
\section{Kubota's formula}
Cauchy's formula expresses the surface area of a convex body $K$ in terms of
the average volume of the projections of $K$ onto codimension $1$ subspaces.  
Kubota's Theorem (see \cite[Chapter 13]{santalo}) is a generalization, which expresses the $k$-th  integral mean
curvature in terms of the average volume of projections of $K$ onto $n-k-1$
dimensional subspaces:
\begin{equation}
\label{kubota}
M_{k}(\partial K) = \dfrac{(n-r-1)\omega_{n-1}}{\omega_{n-k-2}} \fint_{G(n,
    n-k-1)} \vol_{n-k-1}(P_x K) d x,
\end{equation}
where $G(n, n-k-1)$ is the Grassmannian of $n-k-1$ dimensional linear
subspaces of $\mathbb{E}^n,$ and $P_x$ is the projection onto the subspace
$x.$ 

In the special case where $K$ is an ellipsoid $E$ with axes $a_1, \dotsc, a_n,$
Theorem \ref{measproj} gives us several explicit expressions for the integrand
in Kubota's formula. For the purposes of the next Theorem, 
Eq. \eqref{form2} fill be the most useful.
\begin{theorem}
\label{theestimate}
Let $E$ be an ellipsoid with axes $a_1, \dotsc, a_n.$ Let $a_{\mindex{i}},$
for a multindex $\mindex{i} = (i_1, \dotsc, i_{n-k-1})$ be defined as:
\[
a_{\mindex{i}} = \prod_{l=1}^{n-k-1} a_{i_l},
\]
and let 
\[
\mathcal{A} = \sqrt{\sum_{\mindex{i}} a_{\mindex{i}}^2},
\]
where the sum us taken over increasing multindices.
Then 
\begin{equation}
\label{theestineq}
M_k^{(n)}(B^k(1)) \mathcal{A} \leq M_k(E) \leq
\dfrac{M_k(B^n(1))}{\sqrt{\binom{n}{k+1}}} \mathcal{A}.
\end{equation}
\end{theorem}
\begin{proof}
The proof is identical to the proof of Lemma \ref{sharplem}, except we use
$a_{\mindex{i}}$ as variables. With the normalization $\mathcal{A} = 1$
(allowed by homogeneity) we see
  that the maximal case corresponds to the ball of such a radius that 
$a_{\mindex{i}} = 1/\sqrt{\binom{n}{n-k-1}},$ for any multindex $\mindex{i},$
  and the minimum corresponds to $a_1 = \dotso = a_{n-k-1} = 1,$ while
$a_{n-k} = \dotso = a_n = 0.$
\end{proof}

The ratio of the right hand side of the inequality \eqref{theestineq} to the
left hand side is the subject of Theorem \ref{meancrat2} and Remark
\ref{quotest}. As commented in the Remark, the ratio is bounded for any fixed
$k,$ and in the worst case (for $k = n/2$), the ratio grows like $n^{1/4}.$ 

\section{Some exterior algebra}
\label{extalg}
Let $V$ be a vector space, and let $A$ be a linear transformation:
\[
A\in \Hom(V, V).
\]
The exterior power $\bigwedge^k V$ is the vector space generated by
multivectors of the form $v_1 \wedge \dotso \wedge v_k,$ and so we define
\[\bigwedge^k A \in \Hom\left(\bigwedge^k V, \bigwedge^k V\right)\]
by
\[\bigwedge^k A(v_1\wedge\dotso\wedge v_k) = Av_1 \wedge \dotso \wedge A
v_k.\]

From now on, we assume that $V$ is an $n$-dimensional Hilbert space.
The vector space $\bigwedge^k V$ has a standard orthonormal basis: all
multivectors of the form $e_{i_1} \wedge \dotso e_{i_k},$ where 
the $e_{i_l}$ are the standard orthonormal basis vectors in $V,$ and
$i_r \neq i_s,$ for $r \neq s.$ For notational convenience, we will henceforth
denote such multi-indices by bold latin letters. In addition, if $\mathbf{i}$
is a $k$-multindex, we define the $n-k$-multindex $\overline{\mathbf{i}}$ by
\[
e_{\mindex{i}} \wedge e_{\mindexc{i}} = e_1 \wedge e_2 \dotso \wedge
e_n.
\]
\begin{remark}
Riemannian geometers would say that $e_{\mindexc{i}}$ is the image of
$e_{\mindex{i}}$ by the  \emph{Hodge $*$} operator.
\end{remark}
\begin{lemma}
\label{hodgeaux}
If $\mindex{j}\neq \mindex{i},$ then
\[
e_{\mindex{j}}\wedge e_{\mindexc{i}} = 0.
\]
\end{lemma}
\begin{proof}
One of the coordinates of $\mindex{j}$ must be the same as one of the coordinates
of $\mindexc{i}.$
\end{proof}
In the sequel, we will use the following easy observation:
\begin{lemma}
\label{hodge}
Let $\mathbf{v} \in \bigwedge^k V,$ and let $\mathbf{i}$ be a $k$-multindex.
Then 
\[
\langle \mathbf{v}, e_{\mathbf{i}}\rangle = 
\dfrac{\mathbf{v} \wedge e_{\overline{\mathbf{i}}}}{e_1 \wedge \dotso \wedge e_n}.
\]
\end{lemma}
\begin{proof}
Expand $\mathbf{v}$ in coordinates; the result follows immediately from 
Lemma \ref{hodgeaux}.
\end{proof}
\begin{lemma}[Binet-Cauchy formula]
\label{binetcauchy}
Let $A,B \in \Hom(V, V).$ Then, 
\[
\left\langle \bigwedge^k(AB) e_{\mindex{j}}, e_{\mindex{k}} \right\rangle = 
\dfrac{1}{k!}
\sum_{\text{all $k$-multindices $\mindex{i}$}}
\left\langle\bigwedge^k A e_{\mindex{j}}, e_{\mindex{i}}\right\rangle
\left\langle\bigwedge^k B e_{\mindex{i}}, e_{\mindex{k}}\right\rangle.
\]
\end{lemma}
\begin{proof}
This is just the usual formula for matrix multiplication applied in the space
$\bigwedge^k V.$
\end{proof}
We can use the results above to give some identities for projections:
\subsection{On projections}
\begin{theorem}
\label{projlemma}
Let $P$ and $Q$ be such that:
\begin{enumerate}
\item $\rank{P}=k.$
\item $\rank{Q}=n-k.$
\item $P+Q = I.$
\end{enumerate}
Then
\begin{equation}
\label{projeq}
\left\langle \bigwedge^k Pe_{\mindex{i}}, e_{\mindex{i}}\right\rangle =
\left\langle \bigwedge^{n-k} Qe_{\mindexc{i}}, e_{\mindexc{i}}\right\rangle
\end{equation}
\end{theorem}
\begin{proof}
\[
\begin{split}
\left\langle \bigwedge^k Pe_{\mindex{i}}, e_{\mindex{i}}\right\rangle &= \\
\dfrac{\bigwedge^k Pe_{\mindex{i}} \wedge e_{\mindexc{i}}}{e_1\wedge e_n} &= \\
\dfrac{\bigwedge^k Pe_{\mindex{i}} \wedge
  \bigwedge^{n-k}(P+Q)e_{\mindexc{i}}}{e_1\wedge e_n} &= \\
\dfrac{\bigwedge^k Pe_{\mindex{i}} \wedge
  \bigwedge^{n-k}(Q)e_{\mindexc{i}}}{e_1\wedge e_n} &= \\
\dfrac{\bigwedge^k (P+Q)e_{\mindex{i}} \wedge
  \bigwedge^{n-k}(Q)e_{\mindexc{i}}}{e_1\wedge e_n} &= \\
\dfrac{\bigwedge^k e_{\mindex{i}} \wedge
  \bigwedge^{n-k}(Q)e_{\mindexc{i}}}{e_1\wedge e_n} &= \\
\left\langle \bigwedge^{n-k} Qe_{\mindexc{i}}, e_{\mindexc{i}}\right\rangle,
\end{split}
\]
where we have used the observation that $\bigwedge^l P = 0,$ whenever $l >
\rank{P}.$ 
\end{proof}
Suppose $W$ is a subspace of $V,$ and let $w_1, \dotsc, w_k$ be an orthonormal
basis of of $W.$ Let $\Omega$ be the matrix whose columns are the vectors
$(w_1, \dots, w_k, 0, \dotsc, 0)$ (padding $\Omega$ by zeros is not really
necessary, but it will make the sequel slightly simpler notationally).
We then have the following:
\begin{lemma}
\label{projform}
Let $P$ be the orthogonal projection onto $W.$ then 
\[P = \Omega \Omega^t.\]
\end{lemma}
\begin{proof}
The proof is by direct computation: we will show that $Q = \Omega \Omega^t$ is 
the sought-after projector. 
First, let $v$ be orthogonal to all of $W.$ Then, it is clear that $\Omega^t v
= 0,$ and so $Q v = 0.$ Now, consider $Q w_i.$ First, $\Omega^t w_i = e_i.$
Now, for any matrix $A,$ $Ae_i$ is the $i$-th column of $A.$ In particular, 
$\Omega e_i = w_i,$ and so $Q w_i = w_i,$ for all $i.$ It follows that $Q$ is
the sought-after projector.
\end{proof}

\begin{corollary}
\label{projsquare}
Let $P$ and $\Omega$ be as above. Then 
\[
\left\langle \bigwedge^k P e_{\mindex{i}}, e_{\mindex{i}} \right\rangle =
\left\langle \bigwedge^k \Omega e_{\mindex{i}}, e_{\mindex{i}} \right\rangle^2
\]
\end{corollary}
\begin{proof}
This is an immediate consequence of Lemma \ref{projform} above and Lemma
\ref{binetcauchy}. 
\end{proof}

The following can be viewed as a generalization of the Pythagorean theorem:
\begin{theorem}[Generalized Pythagorean Theorem]
Let $\Omega$ be as above. Then the sum of squares of $k\times k$ minors of
$\Omega$ equals $1.$
\end{theorem}
\begin{proof}
By examination of the characteristic polynomial of $P,$ the product of the
non-zero eigenvalues of $P$ equals the sum of the principal $k\times k$ minors,
which, by corollary \ref{projsquare} equals the sum of squares of the
$k\times k$ minors of $\Omega.$ However, since $P$ is a projection, its
non-zero eigenvalues are all equal to  $1.$
\end{proof}

\section{How to compute the volume of a projected ellipsoid}
\label{howtovol}
First, consider a generalized ellipsoid $E(A)$ -- the image of the unit ball
in $\reals^n$ under a linear transformation of rank $k.$ We would like to know
the $k$-dimensional volume of $E(A).$ The simplest situation is when 
\begin{equation}
\label{diagcase}
A_{ij} = \begin{cases}
\lambda_i, \quad i=j, i\leq k,\\
0, \quad \text{otherwise}
\end{cases}
\end{equation}
In this case, 
\begin{equation}
\label{singprodaux}
\vol_k (E(A)) = \kappa_k \prod_{i=1}^k \lambda_i.
\end{equation}
The general case is not much different: Any $A$ of rank $k$  can be written as 
$U \Sigma V,$ where $V\in O(n),$ and $U$ is in $O(k) \subset O(k),$ while
$\Sigma$ is the diagonal matrix of type described in Eq. \eqref{diagcase}; the
diagonal entries of $\Sigma$ are the \emph{singular values} of $A,$ which can
be alternately described as the positive square roots of the (nonzero)
eigenvalues of either $A^t A$ or $A A^t.$
Let us state this as a theorem:
\begin{theorem}
\label{singvalform}
Let $E(A)$ be the image of the unit ball in $\reals^n$ under a transformation
$A$ of rank $k.$ Then 
\begin{equation}
\label{singprod}
\vol_k (E(A)) = \kappa_k \prod_{i=1}^k \sigma_i, 
\end{equation}
where $\sigma_1, \dots, \sigma_k$ are the singular values of $A.$
\end{theorem}

Now, we note that for a matrix $M$ of rank $k,$ the product of the non-zero
eigenvalues of $M$ equals the sum of $k\times k$ principal minors of $M$ (this
is immediate by examining the characteristic polynomial of $M.$ Thus 
Eq. \eqref{singprod} can be rewritten as:
\begin{equation}
\label{detform}
\vol_k^2 (E(A)) = \kappa_k \sum_{\text{principal submatrices $M$ of $A^t A$}}
\det M.
\end{equation}
This last form is superior to Eq. \eqref{singprod}, since 
it expresses the square of the volume as a polynomial in the entries of $A.$
We also note that the $k\times k$ principal minor $M_{\mindex{i}}$ of a matrix
$M$ is something we have already seen:
\[
M_{\mindex{i}} = \left\langle \bigwedge^k M e_{\mindex{i}},e_{\mindex{i}}\right \rangle
\]

\subsection{How do we compute the volume of a projection of an ellipsoid?}
Here we consider a special case: we take a non-degenerate ellipsoid
$E(\mathfrak{A} ),$
where, for simplicity, $\mathfrak{A} = \diag{a_1, \dots, a_n},$ and we would like to
compute the volume of the projection of $E(\mathfrak{A})$ onto a
$k$-dimensional subspace $W$ with the associated  projector $P.$ In other
words, we want to compute the volume of $E(P\mathfrak{A}).$ With the notation 
$A = P \mathfrak{A}$ we note that $A^t A = \mathfrak{A} P \mathfrak{A}.$ To
use the formula \eqref{detform} we first note that if $\mindex{i}$ is a
multindex, then we have the following expression for the minors of $A^t A:$
\begin{equation}
\label{minorform}
(\mathfrak{A} P \mathfrak{A})_{\mindex{i}} = a_{\mindex{i}}^2 P_{\mindex{i}},
\end{equation}
where, if $\mindex{i} = (i_1, \dotsc, i_k),$ then 
\[
a_{\mindex{i}} = \prod_{l=1}^k a_{i_l}.
\]

We then have the following:
\begin{theorem}[Measure of ellipsoid projections]
\label{measproj}
Let $E$ be the ellipsoid 
\[
E = \left\{x \in \reals^n \ \rvert \ \sum_{i=1}^n \dfrac{x_i^2}{a_i^2} \leq
1.\right\}
\]
Let $W$ be a subspace of $\reals^n,$ with an orthonormal basis $w_1, \dots,
w_k,$ while $W^\perp$ has the orthonormal basis $w_{k+1}, \dotsc,  w_n.$ Let
$\Omega$ be the $n \times k$ matrix whose rows are the vectors $w_1, \dots,
w_k,$ while $\Omega^\perp$ be the $n \times (n-k)$ matri whose rows are 
$w_{k+1}, \dotsc, w_n.$ Let the projection onto $W$ be denoted by $P,$ while
the projection onto $W^\perp$ be denoted by $P^\perp.$
Then, the $k$-dimensional volume $\vol_k P(E)$ of the projection of $E$ onto
$W$ can be expressed in any one of the following ways (all the sums below are
taken over \emph{nondecreasing} multindices $\index{i} = (i_1,
\dotsc, i_k),\quad i_1 \leq i_2 \leq \dotso \leq i_k$):
\begin{subequations}
\begin{equation}
\label{form1}
\vol_k P(E) = \kappa_k \sqrt{\sum_{\text{nondecreasing $k$-multindices $\mindex{i}$}} P_{\mindex{i}} a_{\mindex{i}}^2}.
\end{equation}
\begin{equation}
\label{form2}
\vol_k P(E) = \kappa_k \sqrt{\sum_{\text{nondecreasing $k$-multindices
  $\mindex{i}$}} \Omega_{\mindex{i}}^2   a_{\mindex{i}}^2}.
\end{equation}
\begin{equation}
\label{form3}
\vol_k P(E) = \kappa_k \prod_{i=1}^n a_i \sqrt{\sum_{\text{nondecreasing $n-k$-multindices
  $\mindex{i}$}} \dfrac{P^\perp_{\mindex{i}}}{a_{\mindex{i}}^2}}.
\end{equation}
\begin{equation}
\label{form4}
\vol_k P(E) = \omega_k \prod_{i=1}^n a_i \sqrt{\sum_{\text{nondecreasing $n-k$-multindices
  $\mindex{i}$}} \dfrac{(\Omega^\perp_{\mindex{i}})^2}{a_{\mindex{i}}^2}}.
\end{equation}
\end{subequations}
\end{theorem}
\begin{proof}
The expression \eqref{form1} follows immediately from Eq. \eqref{detform}. The
expression \eqref{form2} follows for Eq. \eqref{form1} and Corollary
\ref{projsquare}. The expression \eqref{form3} follows from Eq. \eqref{form1}
and Theorem \ref{projeq}. The expression \eqref{form4} follows from
Eq. \eqref{form3} and Corollary \ref{projsquare}.
\end{proof}
\begin{remark}
The last two expressions (\eqref{form3} and \eqref{form4} in the above theorem
are more useful when $k > n/2;$ the forms \eqref{form2} and \eqref{form4}
are useful when the subspaces are given by their generating vectors, while
\eqref{form1} and \eqref{form3} are more useful when the subspaces are given
by their projectors.
\end{remark}
\begin{example}
\label{dim1}
Suppose $k = 1,$ so we are projecting on a subspace spanned by a (unit) vector 
$v = (v_1, \dotsc, v_n).$ Then, the length of the projection of our ellipsoid is 
(according to Eq. \eqref{form2} is 
\[
\sum_{i=1}^n v_i^2 a_i^2.
\]
\end{example}

\begin{example}
\label{dimn1}
Suppose $k=n-1,$ so we are projecting onto the orthogonal complement of the
subspace subspace spanned by $v$ of the previous example.
Then, the $n-1$-dimensional volume of the projection of our ellipsoid is
(according to Eq. \eqref{form4}):
\[
\kappa_{n-1} \prod_{i=1}^n a_i \sum_{j=1}^n \dfrac{v_j^2}{a_i^2}.
\]
This formula was previously obtained (by completely different methods) by
Connelly and Ostro in \cite{connos}.
\end{example}

\section{Isoperimetric questions}
\label{isoperimetric}
Theorem \ref{theestimate} can be expressed as follows:
\begin{citation}
Let $\mathcal{E}$ be the set of ellipsoids such that the squares
of the $n-k-1$ dimensional volumes of the projections onto coordinate $n-k-1$
dimensional subspaces equals $1$ Then the largest value of
$M_k(E)$ for $E \in \mathcal{E}(\mathcal{A})$ is achieved by the
$n$-dimensional ball (of radius $\binom{n}{n-k-1}^{-1/(n-k-1)},$)
while the minimal value of $M_k(E)$ is achieved by any $n-k-1$ dimensional
ellipsoid parallel to one of the coordinate subspaces.
\end{citation}

It is natural to ask whether the above statement holds with the 
word ``ellipsoid'' replaced by the word ``convex body'' throughout. 
I believe that the answer is in the affirmative, but it is clear that the
methods of this paper do not apply to this question in this generality.

\bibliographystyle{amsplain}

\end{document}